# A BRIEF HISTORY OF
# INFORMATION-BASED COMPLEXITY

Joseph F. Traub
Computer Science Department
Columbia University

This talk is titled a brief history of information-based complexity. Like all memoirs it reflects my personal view of the history of the field. Furthermore, it has a Carnegie Mellon and Columbia University slant. If one of my colleagues were to write a history I'm sure it would differ from this one. I do hope it captures the essence of the field and I apologize to other researchers for any omissions.

I'll begin my history in 1972 when I was Head of the Computer Science Department at Carnegie Mellon University. I received a registered package containing a paper and a letter from someone named Henryk Woźniakowski in Warsaw. I don't recall the date because I didn't realize it was to be the beginning of a life transforming relation which has already lasted 36 years.

The paper, which was titled "Maximal Stationary Iterative Methods for the Solution of Operator Equations" proved conjectures I had framed in the early 60's with a very important difference. My conjectures had been for scalar problems; Woźniakowski proved them for finite-dimensional and infinite-dimensional operator equations.

In a flashback to 1959 I'll tell you why I was so excited by this paper. I had just received my PhD and was working at the Mathematics Research Center at Bell Laboratories. One day a colleague named Joseph Kruskal asked me for advice on how to numerically approximate a zero of a function involving an integral. Since the integral had to be approximated at each iterative step the function was expensive to compute. I could think of a number of ways to solve this problem. What was the optimal algorithm, that is, the method which would minimize the required computational resources? To my surprise there was no theory of optimal algorithms.

Indeed, the phrase computational complexity, which is the study of the minimal computational resources required to solve problems, was not introduced until 1965 by Hartmanis and Stearns [1].

I became fascinated with creating what might be called optimal iteration theory. The initial problem was to solve the scalar nonlinear equation $f(x) = 0$. Assume for simplicity that the zero is simple.

The key insight was that the information used by an iteration determines the maximal order and the most effective methods are iterations of maximal order. The focus was on the information and not the particular algorithm that used the information. Let me give you some concrete examples. You're all familiar with Newton iteration. As you know, it requires the evaluation of f and its first derivative at each step and its order of convergence is two. The order of convergence is a measure of how fast the iteration converges. We can obtain an iteration of order 3 which uses evaluations of f and its first and second derivatives. Generally, there are known methods due to Euler and Chebyshev which use the first s-1 derivatives and are of order s [2, p.81]. These are one-point iterations; that is, all the evaluations are at one point. The question that interested me was whether we could do better. Was it possible for there to exist a one-point iteration of order s that did not require the evaluation of f and its first s-1 derivatives? The answer is no. I called this the fundamental theorem of one-point iterations [2, p.97]. Any one-point iteration of order s must evaluate f and its first s - 1 derivatives.

I don't have to look at the structure of the iteration, only at the information it uses. To put it another way, the maximal order of any one-point iteration using s-1 derivatives is s. For later research on maximal order, see for example [3].

This is a huge simplification. The maximal order of convergence is determined completely by the information available to the iteration, not by its particular form. The significance of iterations of maximal order is that if the cost of information, that is the evaluation of f and its derivatives, is sufficiently large one can neglect the cost of combining the information and the best methods are iterations of maximal order.

You're all familiar with the secant method. Where does it fit in? From an information point of view the secant method evaluates f and reuses one previous evaluation of f. Its order is the golden mean which, as you know, is about 1.62. A method that uses a new evaluation of f and no previous evaluations is of order 1. So the previous evaluation of the secant method adds .62 to the order. There's a method which reuses 2 previous evaluations of f and is of order about 1.84. There's a method which reuses three previous values whose order is about 1.92. These iterations are examples of one point iterations with memory. This and other data suggest that with any finite number of previous values the order will be less than two. That is, the previous evaluations add less than one to the order.

Iterations that use new values of f are a special case. It's natural to consider iterations which use new values of f and its first s - 1 derivatives at a point and reuse any number of previous values. I defined the class of interpolatory iterations and proved that all the old information adds less than one to the order. This is a theorem for interpolatory iteration. I conjectured that this was true for any one-point iteration with memory. This was one of the topics covered in a 1964 monograph called "Iterative Methods for the Solution of Equations" [2]. I'm pleased that it's been reissued by the American Mathematical Society and is still in print.

That's the end of the flashback and I want to return to 1972 when I received the paper from Henryk Woźniakowski. As I mentioned earlier, he attacked the problem of maximal order for finite-dimensional and infinite-dimensional operator equations [4]. He proved the maximal order of interpolatory iteration in the scalar case thus settling the conjecture about one-point iterations with memory. He also proved that in the operator case any one-point iteration of order s requires the evaluation of the first s-1 derivatives.

I invited Woźniakowski to give a talk at a May 1973 Carnegie-Mellon University Symposium. He could not obtain a passport in time to participate. He finally arrived on October 16, 1973. He told me later what led to his paper. He was attending a summer school in Gdansk. Stefan Paszkowski of the University of Wroclaw asked if he'd read my 1964 monograph. He had not but he obtained a copy—the rest is history. Henryk continued to visit me and together with my former PhD student, H.T. Kung, now a chaired professor at Harvard, we continued to work on optimal iteration theory. Then in 1976 there came an event that changed the course of our research.

A PhD student named Arthur Werschulz, now a professor at Fordham University and part of our research group at Columbia, gave a seminar where he used some of the techniques from nonlinear equations to attack the complexity of integration. Our reaction was that integration is inherently different from solving nonlinear equations; one doesn't solve integration iteratively. Because these problems are so different there must be a very general structure which underlies this and many other problems. Henryk and I always maintained long lists of research ideas. But we were so interested in this issue that we called it the S problem which stood for Special problem.

Our search for the general structure led to the 1980 monograph "A General Theory of Optimal Algorithms" [5]. We developed the theory over normed linear spaces with applications to problems such as approximation and linear partial differential equations. We confined ourselves to the worst case setting. That is we guaranteed an approximation for all inputs in a class.

We called the field analytic complexity. This was to differentiate it from algebraic complexity which was a very active research area in the late 60s and 70s. Algebraic complexity deals with problems such as the complexity of matrix multiplication which can be solved exactly while analytic complexity deals primarily with problems from analysis which cannot be solved exactly.

Part B of our 1980 book deals with an iterative information model. It turns out that this material is conceptually and technically more difficult. It was a historical accident which I've told you about earlier that we started with the study of nonlinear equations.

We also gave a brief history of the precursors to the general theory. I'd like to mention a few of the earlier results. These all dealt with specific problems and did not attempt a general theory. The earliest paper which we discovered only recently is by Richard von Mises [6] which was published in 1933 in the Zeitschrift für Angewandte Mathematik und Mechanik. He considered univariate integration with fixed nodes and found the best weights. Arthur Sard authored a series of papers starting in 1949 and a monograph [7] in which he studied optimal algorithms for univariate quadrature with fixed nodes. He discussed extending his results to the approximation of linear functionals. Sard was apparently not familiar with the paper of von Mises. In 1950, Sergei Nikolskij [8] independently studied univariate integration but permitted the evaluation points to be optimally chosen. Another 1950 paper on univariate integration was written by Hans Bueckner [9]. In a series of remarkable papers starting in 1959 Nikolaj Bakhvalov [10] studied optimal methods for multivariate integrals and obtained lower bounds on the error.

All these authors assume linear algorithms; that is, algorithms that are a linear combination of the information. Then in 1965 Sergei Smolyak [11] proved that for convex and balanced sets the optimal algorithm for the approximation of linear functionals is linear. Therefore, the assumption of linear algorithm is often not needed.

As you know the optimal strategy for approximating a zero of a continuous scalar function with a sign change is bisection. What about approximating a maximum of a unimodal function, that is a function which has only one maximum. In a 1953 publication Jack Kiefer [12] proved that if function evaluations are used then Fibonacci search is optimal. This was his 1948 MIT Master's thesis which was only published later with the encouragement of Jacob Wolfowitz. The previous work on optimal algorithms was for linear problems such as integration and approximation. To the best of my knowledge this was the first result for a nonlinear problem.

In 1983 Grzegorz Wasilkowski joined Henryk and me to write the monograph "Complexity and Information" [13]. We showed that uncertainty could be measured without a norm or metric. We decided to rename the field ε-complexity.

One day my wife, Pamela McCorduck, asked me why ε-complexity. I explained that denotes a small quantity. She did not seem impressed. Since Pamela is the author of numerous books I took her lack of enthusiasm seriously and started thinking about a new name. One day I was chatting with my friend, Richard Karp, who as you know was a pioneer in the study of NP-completeness. I mentioned to Richard that key concepts were information and complexity and he suggested information-based complexity which we adopted as the name of the field. For brevity I will often refer to the field as IBC.

The Journal of Complexity was born in 1985. To the best of my knowledge it was the first journal with complexity as its title. In preparing for writing this talk I took a look at Volume 1. There were thirteen people on the Editorial Board. They included three Nobel Laureates (Kenneth Arrow, Gerard Debreu, and Leonid Hurwicz, who is the most recent winner of the Economics Prize), one Fields medalist (Steven Smale), two Turing Prize winners (Michael Rabin and Richard Karp), the founder of Mathematica (Stephen Wolfram), one of the pioneers of algebraic complexity (Shmuel Winograd), the current President of Tel Aviv University (Zvi Galil), a chaired Harvard professor (H.T.Kung) and the recipient of an honorary doctorate from Friedrich Schiller University Jena (Henryk Woźniakowski).

The first volume consisted of two issues containing 285 pages. All the papers were from the Symposium on the Complexity of Approximately Solved Problems held at Columbia in April 1985. Jumping forward to the present the Journal of Complexity now publishes some 1000 pages annually in six issues.

In 1988 Erich Novak published "Deterministic and Stochastic Error Bounds in Numerical Analysis" [14] based on his Habilitation thesis. He studies worst case error bounds which he connects with Kolmogorov n-widths. He also studies error bounds in the randomized and average case settings. The theory is applied to problems such as approximation, optimization, and integration.

In 1988 Wasilkowski, Woźniakowski and I published "Information-Based Complexity" [15]. It integrates the work of numerous researchers and reports many new results. The theory is developed over abstract linear spaces, usually Hilbert or Banach spaces. The worst, average, probabilistic, and asymptotic settings are analyzed. Numerous applications are also presented; these are developed over function spaces. Applications include function approximation, linear partial differential equations, integral equations, ordinary differential equations, large linear systems, and ill-posed problems.

Information-based complexity is defined as the branch of computational complexity that deals with the intrinsic difficulty of the approximate solution of problems for which information is partial, contaminated, and priced. To motivate this characterization consider the numerical solution of a partial differential equation. The coefficients and the initial or boundary values are specified by functions. Since functions cannot be input to a digital computer we have to discretize them by, for example, evaluating them at a finite number of points. Thus a function is represented by a vector of numbers in the computer. There are usually an infinite number of functions which are all represented by the same vector; the mapping is many to one. We say the information about the mathematical input is partial. In addition, there will be round-off errors in evaluating the function and so the information is contaminated. If information is partial and contaminated the problem can only be approximately solved. Finally we'll be charged for evaluating the functions. So the information is priced. Indeed for many problems the cost of the information dominates the cost of combining the information to get an answer.

The next decade was one of rapid progress in IBC which I'll indicate by briefly summarizing five monographs published during that period.

In 1991 Arthur Werschulz published "The Computational Complexity of Differential and Integral Equations: An Information-Based Approach" [16]. Werschulz studies algorithms and complexity of elliptic partial differential equations in the worst case setting. He also studies Fredholm integral equations of the second kind as well as ill-posed problems. In addition, there's a chapter on the average case setting.

In 1996 Leszek Plaskota published "Noisy Information and Computational Complexity" [17]. Plaskota studies both bounded and stochastic noise. Before his work the study of noisy information had lagged due, at least in part, to the technical difficulties.

1998 saw the publication of "Complexity and Information" [18] by Traub and Werschulz. This monograph is a greatly expanded and updated version of a series of lectures I gave in 1993 in Pisa at the invitation of the Accademia Nazionale dei Lincee. It starts with an introduction to IBC and then moves to a variety of topics including very high-dimensional integration and mathematical finance, complexity of path integration, and assigning values to mathematical hypotheses. It concludes with a bibliography of over 400 papers and books published since 1987.

Klaus Ritter's monograph on "Average-Case Analysis of Numerical Problems" [19], which was based on his Habilitation thesis, appeared in 2000. The book provides a survey of results that were mainly obtained in the last ten years as well as many new results. Background material on reproducing kernel Hilbert spaces, random fields, and measures on function spaces is included.

2001 saw the publication of "Optimal Solution of Nonlinear Equations" [20] by Kris Sikorski. The monograph studies algorithms and complexity in the worst case setting. Topics include nonlinear equations, fixed points of contractive and noncontractive mappings, and topological degree.

I'll now return to 1991 when we held the first Schloss Dagstuhl Seminar on Algorithms and Complexity for Continuous Problems. The Schloss Dagstuhl Seminars are the computer science equivalent to the Oberwohlfach meetings in mathematics. In 2009 we will hold our ninth Seminar which may be a record. The Seminars are limited to 40 participants. As the field has grown it has become increasingly difficult for the Organizing Committees to issue only enough invitations to have some 40 participants.

In 1991 Woźniakowski published a paper "Average Case Complexity of Multivariate Integration" [21] in the Bulletin of the American Mathematical Society which was to be a rich source of new directions for IBC. He showed that the optimal points were related to the low-discrepancy points which had been extensively studied by number theorists including Fields Medalist Klaus Roth.

Quasi-Monte Carlo methods are deterministic methods based on low discrepancy points. A Columbia student, Spassimir Paskov, found empirically that Quasi-Monte Carlo was far superior to Monte Carlo for real-world problems in computational finance [22]. In trying to understand why, Ian Sloan and Henryk Woźniakowski introduced the idea of weighted spaces [23]. Discrepancy theory itself has also been a rich source of IBC problems.

1994 saw the first MCQMC Conference organized by Harald Niederreiter. This biennial conference is devoted to Monte Carlo and Quasi-Monte Carlo methods. Many IBC topics and researchers are represented in MCQMC. The 2008 conference will be in Montreal with Stefan Heinrich as Chair of the Steering Committee. The 2010 conference will be in Warsaw with Henryk as the Chair.

1995 saw the first Conference on the Foundations of Computational Mathematics. These triennial meetings always feature an IBC Workshop and an IBC plenary speaker. The 2008 conference will be in Hong Kong and the IBC plenary speaker will be Henryk.

1996 saw the creation of the Best Paper Award of the Journal of Complexity which carries a $3000 prize and a plaque. Roughly half the winning papers have been in IBC. Incidentally, Erich Novak won in 2001 for a paper whose significance I'll describe later. The only double winner is Stefan Heinrich, University of Kaiserslautern. Heinrich was a co-winner in 1998 and won again in 2004.

Since I'm on this topic I'll mention two more IBC Prizes. The Prize for Achievement in Information-Based Complexity was created in 1999. The winners of this annual prize have been Erich Novak; Sergei Pereverzev, Johann Radon Institute for Computational and Applied Mathematics, Austrian Academy of Science; Grzegorz Wasilkowski, University of Kentucky; Stefan Heinrich, University of Kaiserslautern; Arthur Werschulz, Fordham University; Peter Mathe, Weierstrass Institute for Applied Analysis and Stochastics, Berlin; Ian Sloan, University of New South Wales; Leszek Plaskota, University of Warsaw; Klaus Ritter, TU Darmstadt; and Anargyros Papageorgiou, Columbia University.

In 2003 we created a third annual prize, the Information-Based Complexity Young Researcher Award for researchers who have not yet reached their 35th birthday [24]. The recipients to date have been Frances Kuo, University of New South Wales; Christiane Lemieux, University of Calgary; Josef Dick, University of New SouthWales; Friedrich Pillichshammer, University of Linz; Jakob Creutzig, TU Darmstadt; Dirk Nuyens, Catholic University, Leuven; and Andreas Neuenkirch, University of Frankfurt.

Another new direction for IBC was initiated by Erich Novak [25] with the publication in 2001 of "Quantum Complexity of Integration". Until this seminal paper, only the complexity of discrete problems on quantum computers had been studied. Novak studied multivariate integration over Hölder classes. He proved exponential quantum speedups over the classical worst case and quadratic speedups over the classical randomized case.

Since then there has been much progress on IBC problems in the quantum setting with many surprising results mostly due to Stefan Heinrich (see for example [26]) and also due to Papageorgiou and Woźniakowski (see for example [27]).

This talk is titled "A Brief History of Information-Based Complexity". I want to end with the present and the future. The European Mathematical Society recently published Volume I of the "Tractability of Multivariate Problems" [28] by Erich Novak and Henryk Woźniakowski. This volume is some 400 pages in length. Volume II, of similar length, is in progress. There is a huge literature on the complexity of problems in d variables. The complexity bounds are usually sharp with respect to $1/\varepsilon$ where $\varepsilon$ is the error threshold but have, unfortunately, unknown dependence on d. To determine if a problem is tractable the dependence on both $1/\varepsilon$ and d must be determined. Tractability requires new proof techniques to obtain sharp bounds on d. There are many surprising results.

Volume I lists 30 open problems which continues an IBC tradition. Dozens of open questions have been listed in many IBC papers and books. Almost all are still open. Many years ago I gave a talk at MIT which I concluded with half a dozen open questions. Marvin Minsky was in the audience and told me he always saves open questions for his students. I told him there were many more open questions where these came from. Why are there so many open questions? IBC is a relatively young field that covers a huge area of optimal algorithms and complexity for continuous mathematics.

What are some of the future directions? I believe the three volume monograph by Novak and Woźniakowski opens up a whole new area of investigation.

Another huge area for research is that of problems specified by nonlinear operators. To date, much of IBC deals with linear operators and their applications such as integration, approximation, integral equations and linear partial differential equations. Attacking problems defined by nonlinear operators will present us with entirely new challenges.

We've come a long way starting with specific problems such as univariate integration and the solution of scalar nonlinear equations and progressing to a general abstract theory with applications ranging from discrepancy theory and computational finance to quantum computing. I believe the next 50 years will see even greater progress.